\documentclass[12pt]{amsart}

\usepackage{amsfonts}
\usepackage{amssymb,amsmath,amsthm}
\usepackage{verbatim}

\theoremstyle{definition}

\newtheorem{defi}{Definition}[section]
\theoremstyle{theorem}

\newtheorem{thm}{Theorem}[section]
\newtheorem{prop}{Proposition}[section]
\newtheorem{cor}{Corollary}[section]
\newtheorem{lm}{Lemma}[section]
\theoremstyle{remark}

\newtheorem{rem}{Remark}
\newtheorem{ex}{Example}

\long\def\symbolfootnote[#1]#2{\begingroup%
\def\thefootnote{\fnsymbol{footnote}}\footnote[#1]{#2}\endgroup}

\begin{document}
\title{Quaternionic-contact hypersurfaces}
\author{David Duchemin}

\maketitle

\symbolfootnote[0]{The author is partially supported by a CIRGET
    fellowship and by the Chaire  de recherche du Canada en alg\`ebre,
combinatoire, et informatique math\'ematique de l'UQAM.}

\begin{abstract}
In this text, we prove that every quaternionic-contact
structure can be embedded in a quaternionic manifold.
\end{abstract}

\section{Introduction}
There has been a great deal of interest recently in the study of
special classes of complete Einstein metrics whose behavior at
infinity  looks like the hyperbolic $\mathbb{K}$-space with
$\mathbb{K} = \mathbb{R}$, $\mathbb{C}$ or $\mathbb{H}$. Among these,
one finds some negatively curved K\"ahler-Einstein metrics constructed
by Cheng and Yau (\cite{Che80}) on bounded strictly pseudo-convex domains of $\mathbb{C}^n$ and
whose conformal infinity is a strictly pseudo-convexe
CR-manifold with the strictly pseudo-convexe CR-structure induced by
the ambiant complex structure of $\mathbb{C}^n$: one says that this
CR-structure is realizable.

 The standard example of a such metric is the complex hyperbolic metric
$g_{\mathbb{C}}$ on the unit ball of $\mathbb{C}^n$ that is explicitely
 given by
$$g_{\mathbb{C}} = \frac{ 4 euc }{\rho} + \frac{1}{\rho^2} ( d \rho^2
+ (I d \rho )^2 ) $$
with $euc$ being the euclidean metric of $\mathbb{C}^n$, $\rho = 1
- |x|^2$ and $I$ being the complex structure of $\mathbb{C}^n$. The
conformal infinity is the standard CR-structure of $S^{2n+1}$ with
contact distribution the maximal $I$-invariant subspace $H = \ker I d
\rho$ of $TS^{2n+1}$ and where the almost complex structure
on $H$ is the restriction of $I$.

 One knows that all
strictly pseudoconvex CR-structures of dimensions at least 7 are locally
realizable in $\mathbb{C}^n$ (\cite{Kur82}, \cite{Aka87}) but there are strictly pseudo-convexe CR 3-manifolds that are
not realizable, even locally (\cite{Nir73}).

 In this note, we are interested in a similar problem concerning the conformal infinity of metrics modelled
on the quaternionic-hyperbolic metric
$g_{\mathbb{H}}$ with Levi-Civita connection $\nabla^{\mathbb{H}}$
  and defined on the $4(n+1)$-ball of
$\mathbb{H}^{n+1}$ by the formula
$$g_{\mathbb{H}} = \frac{4 euc}{\rho} + \frac{1}{\rho^2} ( d\rho^2 + (
I_1   d \rho )^2 + ( I_2 d \rho )^2 + ( I_3 d \rho )^2 ) \, , $$
where $euc$ is the eulidean metric on $\mathbb{H}^{n+1} \simeq
\mathbb{R}^{4(n+1)} $, $\rho =1 - |x|^2$ and $I_1, I_2, I_3$ are the endomorphisms
obtained by right-multiplication by the purely imaginary quaternions
$i,j$ and $k$. The metric $\rho^2 g_{\mathbb{H}}|_ {T S^{4n+3}}$ is
degenerate, and its kernel $H^{can} = \cap_{i=1}^3 \ker (I_i d\rho) $ satisfies 
$$ d( I_i d \rho ) |_{H^{can}} =
4 euc|_{H^{can}} ( I_i \cdot, \cdot ) \,  $$
and defines what we call a quaternionic contact distribution:

\begin{defi}[\cite{Biq00}] Let $M$ be a smooth manifold of dimension $4n+3$. A codimension 3
distribution $H$ on $M$ is quaternionic contact (QC) if there exists
a metric $g_H$ on $H$ such that one can find locally defined 1-forms
$\eta_1$, $\eta_2$ and $\eta_3$ vanishing on $H$ and an induced pointwise quaternionic structure $(I_1,
I_2,I_3)$ on $H$ (i.e. $I_i^2 = -1$ and $I_1 I_2 = - I_2 I_1 = I_3$ ) with
$$d \eta_i |_H  = g_H (I_i \cdot, \cdot ) \, .  $$ 
\end{defi}
In this case, the conformal class of $g_H$ is totally determined by the distribution
$H$. Remark that in \cite{Ale05}, Alekseevsky and Kamishima introduced a notion
of quaternionic CR-structure that can be defined as a QC-distribution
with an induced quaternionic structure $(I_1, I_2, I_3)$ of integable
almost-complex structures.

Let $N$ be a smooth manifold with boundary $M$ admitting a quaternionic
contact distribution $H$. A metric $g$ defined on a neighboorhoud $ M
\times ] 0, a [ \subset N $ of $M$ with coordinates $(x,\rho)$ is
    called asymptotically hyperbolic
    quaternionic with conformal infinity $H$ if 
$$g \sim \frac{1}{\rho} g_H + \frac{1}{\rho^2}( d\rho^2 + \eta_1^2 + \eta^2_2 +
      \eta_3^2) $$
when $\rho$ goes to zero. One can prove that if the dimension of $M$
      is greater than 7, then every QC-distribution $H$ is the conformal infinity of an
      asymptotically hyperbolic
      quaternionic-K\"ahler metric (AHQK metric), \cite{Biq00}. If $dim(M) = 7$, the author found
      a criterion for a given QC-distribution be the conformal infinity of a
quaternionic-K\"ahler manifold, \cite{Duc04}. This condition
      corresponds to the existence of a CR-integrable twistor space
      and one says that the QC-distribution is integrable in this case.

The space $\mathbb{H}^{n+1}$ is endowed with a quaternionic structure,
i.e. a $GL(n+1, \mathbb{H})Sp(1)$-structure with a torsion-free
connection which can be chosen to be the canonical flat connection $\nabla^0$ in
the case of $\mathbb{H}^n$. The distribution $H^{can}$ is the maximal $\{ I_i
\}_{i=1,2,3}$-invariant subspace of $TS^{4n+3}$ and the pointwise quaternionic structure of $H^{can}$ is the restriction of the ambient
pointwise quaternionic structure of $\mathbb{H}^{n+1}$, one says that
$S^{4n+3}$ is a quaternionic contact hypersurface of $H^{n+1}$.

\begin{defi} Let $M$ be a real hypersurface in a
quaternionic manifold $(N, \mathcal{Q})$ and let $H$ be the maximal
$\mathcal{Q}$-invariant subspace of $TN$. The hypersurface $M$ is called a
QC-hypersurface of $N$ if $H$ is quaternionic-contact with
induced pointwise quaternionic structure that coincides with the restriction to
$H$ of
the elements of $\mathcal{Q}$.

\end{defi}

The aim of this note is to investigate to what extent a given QC
manifold can be realized as a real hypersurface of a quaternionic
manifold and thus extending the canonical example of the realization
of $H^{can}$ in $\mathbb{H}^{n+1}$. In fact, the analogy can be made
more precise at the level of the connections. Indeed, the connection
$\nabla^{\mathbb{H}}$ is quaternionic on the unit ball with a pole of
  order $1$ along $S^{4n+3}$. The difference $ \nabla^{\mathbb{H}} -
  \nabla^0$ is a tensor $a^{-\frac{d \rho}{\rho}}$ where 
$a : (\mathbb{H}^{n+1} )^* \rightarrow (\mathbb{H}^{n+1}) ^* \otimes End (
  \mathbb{H}^{n+1})$ is linear and depends only of the $GL(n+1,
  \mathbb{H}) Sp(1)$-structure of $\mathbb{H}^{n+1}$, i.e. $a$ can be
  defined on any quaternionic manifold $M$ and goes from $T^* M$ into
  $T^* M \otimes End (TM)$ (see the next section for an explicit description
  of $a$).

Building on results of \cite{Biq00}, we prove the following theorem:

\begin{thm}\label{thm principal}
Let $(M^{4n+3}, H) $ be an integrable QC-manifold, $n\geq 1$. There exists a
quaternionic manifold $(N_0^{4n+4}, \mathcal{Q})$ such that:
\item[(i)] $M$ is a QC-hypersurface of $N_0$.
\item[(ii)] $M$ separates $N$ into two quaternionic manifolds $N^+_0$
and $N^-_0$ such that $N^+_0$ has a definite-positive
quaternionic-K\"ahler metric with conformal infinity $H$, Levi-Civita
connection $\nabla^{\mathbb{H}}$ and $N^-_0$ has
a quaternionic-K\"ahler metric with signature (4, 4n) and conformal
infinity $H$. 
\item[(iii)] If $\rho$ is a defining function of $M$, then
  $\nabla^{\mathbb{H}} - a^{- \frac{d\rho}{\rho}}$ extends to a smooth
  quaternionic connection on $N$. 
\end{thm}

When $n=0$, $M$ is a conformal 3-manifold which is the conformal
infinity of a unique self-dual Einstein metric. The conformal class of
this metric admits a prolongation in such a way that $M$ becomes an
hypersurface in a self-dual conformal 4-manifold. In this setting,
the theorem proved in this paper appears to be a generalization of this
fact.

On the other hand, an hypersurface $M$ in a 4-dimensional conformal
manifold $(N, [g])$ defines the data of a conformal metric and a second
fundamental form that vanish iff there exists an Einstein metric on
$N-M$ in the conformal class of $[g]$ with conformal infinity
$M$. LeBrun proved that with the data of a conformal metric $[h]$ and a second
fundamental form $\Omega$ on a 3-manifold $M$, one can construct an embedding
of $M$ into a self-dual 4-manifold inducing $( [h], \Omega )$ on
$M$. In this paper, we generalize the notion of conformal second
fundamental form to the case of a QC-hypersurface and prove that it
vanishes for the embedding given by Theorem 1.1. In particular, if $M$
is a QC-hypersurface with non-vanishing second fundamental form in a
quaternionic manifold $N$ with boundary $M$, there does not exist any AHQK-metric on
$N-M$, compatible with the quaternionic structure of $N-M$ and with
conformal infinity $M$.

We now describe briefly the organization of the paper. In the first
section, we recall some basic facts about quaternionic manifolds and give
the local description of a QC-hypersurface and thus define the notion
of weakly quaternionic contact manifolds.

 In the next section, we
define the integrability of a QC-distribution and show that the
QC-distribution of a QC-hypersurface is integrable. 

In the third section, we give the definition of the twistor space of a
QC-structure and prove Theorem 1.1.

The aim of the last section is to define the second fundamental forms
of a QC-hypersurface and prove that they vanish under the hypothesis of Theorem 1.1. 

\bigskip

The author gratefully acknowledge many useful conversations with
Vestislav Apostolov during the preparation of this work.

\section{Hypersurfaces in quaternionic manifolds}

In this preliminary section, I describe some basic
facts in quaternionic geometry, see \cite{Bes87}[p.410] and
\cite{Sal86} for more details. I give also the local description of a
quaternionic contact hypersurface.

\subsection{Quaternionic manifolds}

\begin{defi}
An almost quaternionic manifold is a $4n$-dimensional manifold endowed
with a $GL(n, \mathbb{H} )Sp(1) $-structure. A quaternionic manifold is an almost-quaternionic manifold admitting a
torsion-free $GL(n, \mathbb{H})Sp(1)$-connection.
\end{defi}

\begin{rem}
An almost quaternionic manifold $M$ is the data of a
sub-bundle $\mathcal{Q} \subset End (TM)$, locally spanned by a
pointwise quaternionic structure $(I_1, I_2, I_3)$. The manifold $(M, \mathcal{Q})$ is quaternionic if there exists a torsion-free connection that preserves $\mathcal{Q}$. 
\end{rem}

Let $(M, \mathcal{Q})$ be a quaternionic manifold, and let $\nabla$ be a torsion-free
connection preserving $\mathcal{Q}$. If
$\nabla'$ is another torsion free connection that preserves $\mathcal{Q}$, then $\nabla' =
\nabla + a$ where $a$ is in the kernel $\ker \partial $ of the torsion map
$$\partial : \Lambda^1 \otimes ( \mathfrak{gl}(n , \mathbb{H}  )
\oplus \mathfrak{sp}(1) ) \rightarrow \Lambda^2 \otimes \Lambda^1 \,
. $$
Let $\mathfrak{g}$ be the Lie algebra $ \mathfrak{gl}(n , \mathbb{H}  )
\oplus \mathfrak{sp}(1) $.

\begin{lm}[\cite{Sal86}] \label{changedeconnexionquaternionic}
The kernel of the torsion map $\partial$ is the set $\{ a^{\theta} \in
\Lambda^1 \otimes \mathfrak{g}\, , \; \theta \in \Lambda^1\} $, where
$a^{\theta}$ is defined by
$$ a_X^{\theta} Y = \theta (X) Y  + \theta(Y) X - \sum_{i=1}^3 ( \theta(I_iX)
I_i Y + \theta (I_i Y ) I_i X ). $$

\end{lm}

\begin{proof}
It is well known that the restriction
 of $\partial $ to $\Lambda^1 \otimes \mathfrak{sp}(1)$ is
 injective (this follows from the unicity of a
 Levi-Civita connection). Moreover, the only common irreducible summand in the
$GL(n, \mathbb{H} )Sp(1)$ decomposition of $\Lambda^1 \otimes
 \mathfrak{sp}(1) $ and $\Lambda^1 \otimes \mathfrak{gl}(n, \mathbb{H}
 )$ is $\Lambda^1$, where the embedding of $\Lambda^1$ in $ \Lambda^1 \otimes
\mathfrak{sp}(1)$ is given by $\theta \mapsto \overline{\theta} =  \sum_i \theta \circ I_i \otimes
I_i$. The torsion map $\partial$ is $GL(n, \mathbb{H} )
Sp(1)$-equivariante, hence if $\partial a = 0$ the $\Lambda^1
\otimes \mathfrak{sp}(1)$-part $\overline{\theta}$ of $a$ must live in $\Lambda^1$. Now,
if $a'$ is the $ \Lambda^1 \otimes \mathfrak{gl}(n, \mathbb{H} )
$-part of $a$, we obtain
$$\begin{array}{rcl}
a'_X I_3 Y & = & I_1 a'_X I_2 Y = I_1 a'_{I_2 Y} X - I_1  \partial
\overline{\theta} _X I_2 Y \\
& = & a'_{I_1 X} I_2 Y - \partial  \overline
{\theta }_{I_2 Y} I_1 X  - I_1  \partial \overline{\theta}_X  I_2 Y    \\
 & = & I_2 a'_Y I_1 X -I_2 \partial \overline{\theta} _{I_1 X}Y  - \partial  \overline
{\theta }_{I_2 Y} I_1 X  - I_1  \partial \overline{\theta}_X  I_2 Y
 \\
& = & - I_3 a'_X  Y +  I_3 \partial \overline{\theta} _Y  X  -I_2 \partial \overline{\theta} _{I_1 X}Y  - \partial  \overline
{\theta }_{I_2 Y} I_1 X  - I_1  \partial \overline{\theta}_X  I_2 Y
 \, . 
\end{array}$$
Developing the last line of the previous computation gives the lemma.
\end{proof}

If $(M, \mathcal{Q})$ is an almost quaternionic manifold and $\nabla$
is a connection preserving $\mathcal{Q}$, one defines the torsion
$T(\mathcal{Q})$ of
$\mathcal{Q}$ to be the projection of the torsion $T^{\nabla}$ onto
$\Lambda^2 \otimes \Lambda^1 / \partial ( \mathfrak{gl}(n , \mathbb{H}
) \oplus \mathfrak{sp}(1) ) $. It does not depend of the choice of
connection $\nabla$ preserving $\mathcal{Q}$ and it vanishes iff
$\mathcal{Q}$ is quaternionic.

\subsection{The flat model}

In this section, we describe the link between the flat hyperk\"ahler
metric $euc$ on $\mathbb{H}^{n+1} $ and the standard
quaternionic-K\"ahler hyperbolic metric
$g_{\mathbb{H}}$ on the unit ball $B^{4(n+1)} \subset
\mathbb{H}^{n+1}$. Let $\nabla^0$ be the flat (hyperk\"ahler)
connection on $\mathbb{H}^{n+1}$ and $\nabla^{\mathbb{H}}$ be the
Levi-Civita connection of $g_{\mathbb{H}}$. 

\begin{lm}
The connections $\nabla^{\mathbb{H}}$ and $\nabla^0$ are
quaternionic and related by the formula
$$\nabla^{\mathbb{H}} = \nabla^0 + a^{ - \frac{d  \rho }{\rho} } \, .$$
\end{lm}

\begin{proof}Let $\rho = 1 - |x|^2$ and let $\partial_{\rho}$ be the dual vector field of $d \rho$
respectively to $euc$. 

The connections $\nabla^{\mathbb{H}}$ and $\nabla^0$ are quaternionic
hence there exist a 1-form $\theta$ such that $\nabla^{\mathbb{H}}
= \nabla^0 + a^{\theta}$. Let $X \in \mathbb{H}^{n+1}$ satisfying $d
\rho (X) = d \rho (I_1 X) = d \rho (I_2 X) = d \rho (I_3 X) =
0$. Then, the formula $ ( \nabla^0_{\cdot}  g_{\mathbb{H}} +
a^{\theta}_{\cdot} g_{\mathbb{H}} ) (X, X) = \nabla_{\cdot}
^{\mathbb{H}} g_{\mathbb{H}} (X, X) = 0 $ gives
$$ \frac{4}{\rho} euc( X, X) \theta (\cdot ) + \frac{4}{\rho} \theta (X) ( euc ( X,
\cdot ) - \sum_{i=1}^3 euc (I_i X , \cdot ) ) = $$ 
$$-\frac{4}{\rho^2} euc ( X, X ) d \rho
( \cdot )  \,
. $$
 Applying this to $X \neq 0 $ and to $I_i \partial_{\rho}  $ gives
$\theta (X) = 0$ and $\theta ( I_i \partial_\rho )  = 0$ whereas
applying this to $\partial_{\rho}$ gives $\theta ( \partial_\rho ) = -
1/ \rho  $.
\end{proof}

In this description, one sees that $\nabla^{\mathbb{H}}$ admits a
prolongation $\overline{\nabla} = \nabla^0 + a^{ - \frac{d \rho
}{\rho}} $ to
$\mathbb{H}^{n+1}$ with pole along
$S^{4n+3}$ whereas $\nabla^0  = \overline{\nabla} + a^{\frac{d
\rho}{\rho}}$ is smooth on all $\mathbb{H}^{n+1}$. We will keep this
description in mind in order to prove the
theorem \ref{thm principal}. Indeed, we will prove that if $\nabla$ is
the connection of the AHQK-metric with given boundary $(M, H)$ and
$\rho$ is a defining function of $M$,
then both the quaternionic structure and $\nabla + a^{\frac{d \rho}{\rho}}$
have a smooth prolongation into a
neighbourhood of $M$.

\begin{rem}
As a quaternionic manifold, $\mathbb{H}^{n+1}$ compactifies to $\mathbb{H}P^{n+1}$.
\end{rem}

\begin{rem}
There is another quaternionic-K\"ahler metric defined on
$\mathbb{H}^{n+1}$, with positive scalar curvature, which comes from
the embedding $ \mathbb{H}^{n+1} \hookrightarrow \mathbb{H}P^{n+1} $, $x
\mapsto [ 1 , x ] $. This metric can be written as
$$g _+ = \frac{4 \, euc }{1 +|x|^2 } - \frac{1}{(1+|x|^2)^2} ( (d\rho  )^2
+ ( I_1 d \rho )^2 + (I_2 d \rho )^2 + ( I_3 d \rho^2 ) ) \, $$
on $\mathbb{H}^{n+1}$.

\end{rem}

\subsection{Local description}

In this section, we give the fundamental property of a QC-hypersurface
and discuss the general situation of a real hypersurface in a
quaternionic manifold.

\begin{prop}
Let $(N, \mathcal{Q})$ be a quaternionic manifold and $\nabla$ be a
torsion-free connection preserving $\mathcal{Q}$. Suppose $f : M \rightarrow
\mathbb{R}$ is a smooth function with non-vanishing differential $d_x
f$ for all $x \in M= f^{-1 } ( 0 )$. Then $M$ is a QC
hypersurface of $N$ iff $\nabla df $ defines a positive or negative definite metric on the maximal
$\mathcal{Q}$-invariant subspace $H$ of $TM$ and $ \nabla df (I X, I
Y) = \nabla df (X, Y)   $ for all $X, Y \in H$ and $I \in \mathcal{Q}$.
\end{prop}

\begin{proof} Assume first that $M=f^{-1}(0)$ is a QC-hypersurface in
$N$ and that $(I_1, I_2, I_3 )$ is a local choice of quaternionic
structure defined in a
neighbourhood of $p \in M$. The QC distribution is the distribution $H = \cap_i \ker df
\circ I_i $ in $TM$. By hypothesis, there exists a metric $g$ on $H$
and $J_i \in Vect ( I_1, I_2, I_3)$ such that $d (df \circ I_i ) |_U =
g ( J_i \cdot , \cdot )$. The connection $\nabla$ is torsion-free, hence $A = \nabla
df$ is  symmetric. We obtain for $X, Y \in H$
$$A(X, I_i Y) - A( Y, I_i X) + df ( (\nabla_X I_i ) Y) - df ( (
\nabla_Y I_i ) X ) = g ( J_i X, Y ) \, . $$
Because $\nabla$ preserves the quaternionic structure and $ df ( I_i
X) = 0$ for $X \in H$, we get
$$A( X, I_i Y ) - A ( I_i X  ,Y ) = g ( J_i X, Y ) \, . $$
Therefore, $I_i$ and $J_i$ commute. Using now the fact that $J_i \in
vect (I_j)_j$, we obtain the existence of $\lambda_i \in \mathbb{R}$
such that $J_i = \lambda_i I_i$. On one hand, we have
$$A(X, I_3 Y ) - A ( I_3 X, Y ) = \lambda_3 g ( I_3 X, Y ) \, ,$$
and on the other hand,
$$\begin{array}{rcl}
A(X, I_3 Y)  & = &  A( X, I_1 I_2 Y  )   \\
& = &  \lambda_1 g ( I_1 X, I_2 Y ) + A ( I_1 X , I_2 Y ) \\
 & = & \lambda_1 g ( I_3 X, Y ) - \lambda_2 g ( I_3 X, Y )   -A( I_3
 X, Y  )  
\end{array} $$
 and therefore 
$$2A( X, I_3 Y) = ( \lambda_1 - \lambda_2 - \lambda_3 ) g (I_3 X, Y ).
$$
We thus get $ - 2 A( X, Y)  = ( \lambda-1 - \lambda_2 - \lambda_3 ) g
(X, Y)  \, .$
Hence, there exists a scalar $\lambda$ such that
 $ \lambda_1 - \lambda_2 - \lambda_3 = \lambda$, and by
circular permutation $\lambda_2 - \lambda_3 - 
\lambda_1 = \lambda_3 - \lambda_1 - \lambda_2 = \lambda$.
This gives $\lambda_i = - \lambda \neq 0 $ and $2 A(X, Y) = \lambda g
(X, Y) $.
\end{proof}

\begin{rem}
Assume that $(M, H, g_H)$ is a quaternionic hypersurface in a
quaternionic-K\"ahler 8-manifold $(N, \mathcal{Q},g)$ where
$\mathcal{Q}$ is the quaternionic structure and $g$ is the riemannian
metric. In that case, the conformal class $g_H$ is completely
determined by the quaternionic structure of $H$ and is thus equal to the conformal class of $g|_H$. This is not necessarily the
case in dimension greater than $8$.
\end{rem}

For the sake of completeness, we describe now the structure induced on
a general hypersurface in a quaternionic manifold.

\begin{defi}
Let $M^{4n+3} $ be a smooth manifold. A weakly quaternionic contact
structure on $M$ is the data of a codimension 3 distribution $H$ and a $GL(n
, \mathbb{H}) Sp(1)$-structure $\mathcal{Q}\subset End (H)$ on $H$ such that locally there exist 1-forms
$(\eta_1, \eta_2, \eta_3)$ vanishing on $H$ and a $SO(3)$-basis
$(I_1, I_2, I_3)$ of $\mathcal{Q}$ satisfying :
\begin{itemize}
\item[(i)] $d \eta_i (I_i X, I_i Y)
= d \eta_i  (X, Y)$;
\item[(ii)]  the tensor $g=d \eta_1 ( \cdot , I_1 \cdot ) + d \eta_1 (
 I_2 \cdot , I_3 \cdot  ) $ is non-degenerate on $H$;
\item[(iii)] for all $X, Y \in H $, one has the equalities :
$$ g (X, Y)  = d \eta_2 ( X , I_2 Y) + d \eta_2 ( I_3 X, I_1 Y  ) =
 d \eta_3 ( X, I_3 Y ) + d \eta_3 ( I_1 X, I_2 Y )  \, . $$
\end{itemize}
 The tensor $g$ is symmetric. If it is positive definite, we say
that $( H, \mathcal{Q} )$ is a strictly pseudo-convexe weakly quaternionic contact distribution.
\end{defi}

In this case, and contrary to the quaternionic-contact case,
the quaternionic structure $\mathcal{Q}$ on $H$ is not determined by
the distribution $H$. In order to see that, let us describe a simple linear
algebra example. On $\mathbb{H}$, let $(I_i^+)$ (resp. $(I_i^- )$) be the
action of $i$, $j$ and $k$ on right (resp. on left), and define $w_i^+
= euc ( I_i ^+  \cdot, \cdot
)$ and $w_i^-  = euc ( I_i^- \cdot , \cdot )$. We put for $|\lambda| <1 $, $w_i^{\lambda} = \frac{1}{1-\lambda^2} ( w_i ^+
+ \lambda w_i^- ) $. Then, one has $I_i^+ w_i^{\lambda} = w_i ^{\lambda}$
and
$$w_1^{\lambda} ( \cdot, I_1^+  \cdot ) + w_1^{\lambda} ( I_2 ^+\cdot,
I_3^ + \cdot) = \frac{4}{1-\lambda^2} euc \, , $$
with the other relations obtained by circular permutations. But on the
other hand, we have that $ w_i^{\lambda} \wedge w_j^{\lambda} =2
\delta_{ij} \nu $ where $\nu$ is the volume form of $euc$ so that there exists a quaternionic
triple $(I_i^{\lambda})$ and a metric $g_{\lambda}$ not in the
conformal class of $euc$ such that $
w_i^{\lambda} (\cdot , \cdot ) = g_{\lambda} ( I_i^{\lambda} \cdot,
\cdot )$.

The notion of weakly quaternionic-contact distribution is introduced
in order to describe hypersurfaces in quaternionic manifolds. Indeed,
we have :

\begin{prop}
Let $M^{4n+3}$ be an hypersurface in a quaternionic manifold $(N,
\mathcal{Q})$ and let $H$ be the maximal $\mathcal{Q}$-invariant
subspace of $TM$. We denote still by $\mathcal{Q}$ the set of the
restrictions to $H$ of elements of $\mathcal{Q}$. Let $f$ be a
defining function of $M$, let $(I_1, I_2, I_3)$ be a local
$SO(3)$-trivialization of $\mathcal{Q}$ and put 
$$g = \nabla df |_H + \sum_i I_i \nabla df |_H \, .$$
If $g$ is non-degenerate on $H$, then $(H,
\mathcal{Q} )$ is a weakly quaternionic-contact structure on $M$.
\end{prop}

\begin{proof}
Let $f$ be a function defining $M$ locally, and $\eta_i = - df \circ
I_i $ on $M$. Then, one has $d \eta_1 |_H = \nabla df ( \cdot, \cdot)
+ \nabla df ( I_1 \cdot, I_1 \cdot )$ so that it is $I_1$-invariant
and on $H$,
$$d \eta_1 ( \cdot , I_1 \cdot  ) + d \eta_1 ( I_2 \cdot , I_3 \cdot )
= \nabla df (\cdot , \cdot )
+ \sum_{i=1}^3 \nabla df ( I_i \cdot  , I_i \cdot  )  \,  $$
is $I_i$-invariant for all $i$. The other relations are obtained by
cyclic permutation.
\end{proof}

\section{QC geometry}

This section gives the necessary backgroung about QC-distribution. In
particular, we define the integrability of a QC-distribution and prove
that the QC-distribution of a 7-dimensional QC-hypersurface is
integrable. Then we describe the properties of the Biquard connection
that are used in the fifth section to compare the Biquard
connection of a QC-hypersurface with the ambient quaternionic connection.

\subsection{The group $Sp(n)Sp(1)$}

The group $Sp(1)$ of unit quaternions acts on $\mathbb{H}^n$ by right
multiplication and has centralizer $Sp(n) \subset SO(4n)$. One of the
particular feature of the group
$Sp(n)Sp(1)$ is that it arises in Berger'list of
possible holonomy groups of non locally symmetric riemannian
manifolds. In this paper section, we are mainly interested in describing
some features of the
representations of $Sp(n)Sp(1)$. 

Let $V^{(a_1, \dots, a_n)}$ be the irreducible representation of $Sp(n)$
with highest root $(a_1, \dots, a_n)$. If $\sigma \simeq \mathbb{C}^2$
is the standard representation of $Sp(1)$, then the irreducible
representations of $Sp(n)Sp(1)$ are the tensor products $V^{(a_1,
\dots , a_n)} \otimes \sigma^p $ where $a_1 + \cdots a_n + p $ is even
and $\sigma^p$ is the p-symmetric power of $\sigma$; the real
irreducible representations are the real parts $[ V^{(a_1, \dots, a_n
)} \otimes \sigma^p ] $ of the previous ones.

Following Salamon \cite{Sal89}, we put $\lambda^r_s = V^{(a_1, \dots,
a_n)}$ where $s$ of the $a_i$ are equal to $2$, $r-2s$ are equal to
$1$ and the others are zero and we abbreviate $[\lambda^r_s \otimes
\sigma^p] $ in $[ \lambda^r_s \sigma^p ] $.

With this notation, we have $\mathfrak{sp}(1) = [
\sigma^2 ]$, $\mathfrak{sp}(n) = [\lambda^2_1 ]$ and $[\lambda^2_0 ] $
is the set of symmetric traceless $\mathbb{H}$-linear endomorphisms.

Moreover, we have the decompositions
$$\mathfrak{so}(4n) = \mathfrak{sp}(n)\oplus
 \mathfrak{sp}(1) \oplus [\lambda^2_0 \mathfrak{sp}(1)]  $$
and
 $$\mathfrak{gl}(4n, \mathbb{R}) = \mathfrak{sp}(n)\oplus
 \mathfrak{sp}(1) \oplus [\lambda^2_0 \mathfrak{sp}(1)] \oplus [
 \mathfrak{sp}(n) \mathfrak{sp}(1) ] \oplus [\lambda^2_0 ]  \oplus \mathbb{R}  \, $$
where $\lambda^2_0 = 0$ if $n=1$.

\subsection{Integrablity of a QC structure}

Let $(M, H ,g) $ be a QC distribution and $g$ a compatible metric on
$H$ so that one gets a $Sp(n)Sp(1)$-structure on $H$. Let $(\eta_1, \eta_2, 
\eta_3 )$ be a $SO(3)$-trivialization of the set of 1-forms vanishing
on $H$. If $W$ is a complementary to $H$ and $(R_1, R_2, R_3)$ is a
dual basis of $(\eta_i|_W)$, we put $\alpha_{ij} = \iota_{R_i} d
\eta_j |_H$. Remark that we have the natural identification $W \simeq
\mathfrak{sp}(1)$, $R_i \mapsto I_i$ and one can verify that 
$$T^W = \sum_{i,j =1}^3 (\alpha_{ij} + \alpha_{ji} ) \otimes I_i \otimes
I_j \in H^* \otimes \mathfrak{sp}(1) \odot \mathfrak{sp}(1) \simeq
[\lambda^1 \sigma^1 ] \oplus [ \lambda^1 \sigma^3 ] \oplus [ \lambda^1
  \sigma^5 ]$$ 
is a tensor. Changing $W$ corresponds
to changing $T^W$ by an element in $W^* \otimes H \simeq [\sigma^2]
\otimes [ \lambda^1 \sigma ] $, so that one can prove that there exists
a unique complementary $W^g$ of $H$ such that $T^{W^g} \in [ \lambda^1
  \sigma^5 ]$. The decomposition of $H^* \otimes \mathfrak{sp}(1)
  \odot \mathfrak{sp}(1)$ is explicitely given by :
$$\begin{array}{rcl} 
 \left [ \lambda^1 \sigma^1 \right ] &  = &  \left \{ \sum_i r \otimes I_i \otimes
 I_i \, ,  \;
 r \in H^* \right \} \, , \\
 \left [ \lambda^1 \sigma^3 \right ] & = & \{ \sum_{i, j } ( I_i r_j +
 I_j r_i ) \otimes I_i
  \otimes I_j \, , \; r_i \in H^* \, , \; \sum_i I_i r_i =0 \} \, , \\
  \left [ \lambda^1 \sigma^5 \right ] &  = & \{ \sum_{i,j}
 a_{ij}\otimes  I_i \otimes I_j
   \, , \; a_{ij}=a_{ji} \mbox{ and } \sum_j I_i a_{ij} = 0 \} \, . 
\end{array}
$$
\begin{rem}
The vector fields $R_i$ in $W^g$ are called Reeb vector fields of
$(\eta_1, \eta_2, \eta_3 )$.
\end{rem}

\begin{thm}[\cite{Duc04}] \label{Torsion}
Let $(H,g)$ be a quaternionic contact distribution in a manifold $M$ of
dimension $4n+3$. The tensor $T^{W^g}$ is called the vertical torsion
of $H$. It is conformally invariant and vanishes automatically when
$n>1$. If $n=1$ and $T^{W^g} = 0$, one says
that $H$ is integrable.
\end{thm}

The importance of the integrability condition is enhanced by the
following proposition.

\begin{thm}
A QC-hypersurface $M$ of a quaternionic manifold $N$ of dimension $8$
is integrable.
\end{thm}

\begin{proof}
Assume that $f$ is a defining function for $M$. There exists a vector field $\xi$
defined up to a vector field in $H$ and such
that $df (\xi ) = 1$ and $df (I_i \xi )  = 0$. Moreover, $\nabla df$
is non-degenerate on $H$, hence we can assume that $\nabla df ( \xi,
X) =0 $ for all $X\in H$.  Let us define $\alpha_{ij}
= - i_{I_i \xi} d (df \circ I_j )  |_H$. We have for $X \in H$,
$$\alpha_{ij} (X) =- \nabla df ( I_i \xi , I_j X) +\nabla df ( X, I_j
I_i \xi ) + df (  ( \nabla_X I_j ) I_i \xi ) $$
and therefore,
$$\alpha_{ij} (X) + \alpha_{ji} ( X) = \alpha_i ( I_j X) + \alpha_j
(I_i X) $$
where $\alpha_i (X) = -\nabla df ( I_i \xi , X)$.
\end{proof}

\begin{rem}
The vector field $\xi$ that appears in the previous proof is called
the normal vector field of $f$ along $M$.
\end{rem}

\subsection{The Biquard connection}
In thi section, I describe the connection of Biquard for a QC
distribution. The results I give come from \cite{Biq00} and
\cite{Duc04}. Let $\tilde{g}$ be the metric equals to $ \sum_i \eta_i^2$ on $W^g$, $g$ on $H$ and such
that $H$ and $W^g$ are orthonormal.

\begin{thm} Let $H$ be a QC distribution, integrable if $n=1$ and let
$g$ be an adapted metric. There exists a unique connection $\nabla$ preserving $W^g$, $H$ and
$\tilde{g}$ and satisfying:
\begin{itemize}
\item[(i)] $\nabla$ preserves the $Sp(n)Sp(1)$ structure on $H$.
\item[(ii)] if $R \in W$ and $X \in H$, then the torsion $T(R, X) $ is
in $H$ and $X \mapsto T(R, X) \in ( \mathfrak{sp}(n) \oplus
\mathfrak{sp}(1))^{\perp}$.
\item[(iii)] If $X, Y \in H$, then $T(X, Y) \in W$ and if $R, R' \in
W$ then $T(R, R') \in H$.
\end{itemize}

\end{thm}

\begin{rem}
The point $(ii)$ follows from \ref{Torsion}.  
\end{rem}

\begin{rem}
If $X$ and $Y$ are in $H$, then the torsion $T(X, Y)$ satisfies
$$ T(X, Y) = \sum_i \langle I_i X, Y \rangle  R_i \, .$$
\end{rem}

\begin{rem}
The connection given here differs slightly of the connection of
\cite{Biq00}. In fact, the only differences lies in the terms
$\nabla_R R'$ of the connection, so that the terms $\nabla_X R_i$ and
$\nabla_X I_i$ still coincide when we identify $I_i$ and $R_i$.
\end{rem}

One can prove the following stronger
 result for the torsion:

\begin{prop} Assume that the QC distribution $H$ is integrable if
 $n=1$. The part $T^H$ in $W^* \otimes H^*
\otimes H$ of the torsion of $\nabla$ satisfies $T(R, \cdot ) \in [\lambda^2_0
\mathfrak{sp}(1) ] \oplus [\mathfrak{sp}(n)
\mathfrak{sp}(1)]$ for all $R \in W$ and  
$$ T^H \in [\lambda^2_0] \oplus [ \mathfrak{sp}(n) \mathfrak{sp}(1) ]
\subset \mathfrak{sp}(1) \otimes  (  [\lambda^2_0
\mathfrak{sp}(1) ] \oplus [\mathfrak{sp}(n)
\mathfrak{sp}(1)] ) \, . $$ 
\end{prop}

\begin{rem}
If $n=1$, then $\lambda^2_0 = 0$, so that if $H$ is integrable, then
$T^H \in \mathfrak{sp}(1)\mathfrak{sp}(1)$.
\end{rem}

\begin{rem}
The previous proposition implies the existence of $\tau \in [
\lambda^2_0 ]$ and $\tau_k \in \mathfrak{sp}(n)$ such that
$$T_{R_i}^H = I_i \tau + \sum_{k,j=1}^3 \varepsilon^{ijk} I_j \tau_k \, . $$
\end{rem}

\section{Twistor spaces}

In this section, we will prove theorem 1.1. In a first part, we
recall the definition of the twistor space of a QC-distribution $H$
(integrable in dimension $7$) and
the properties of the twistor space of the AHQK metric which has
conformal infinity $H$.

\subsection{The twistor space of a QC-distribution}\label{section-twist}
Let $(M^{4n+3}, H)$ be a QC-contact distribution that is integrable if
$n=1$, and let $g$ be a compatible metric on $H$. Let $(I_1, I_2,
I_3)$ be is a local quaternionic structure on $H$ and
$(\eta_i)_{i=1,2,3}$ such that $d \eta_i ( \cdot, \cdot)
= g (I_i \cdot, \cdot)$ on $H$. Define
$$\mathcal{T} = \{ x_1 I_1 +x_2 I_2 + x_3 I_3 \, , \; x_1^2 + x_2^2 +
x_3 ^2 = 1 \}  $$
the set of compatible almost complex structures on $H$ and let $\pi : \mathcal{T} 
\rightarrow M$.  The space $\mathcal{T}$ is a 2-sphere bundle
over $M$, and carries a 1-form $ \eta =
x_1 \pi^* \eta_1 + x_2 \pi^* \eta_2 + x_3 \pi^* \eta_3 $ defined up to
a conformal factor. At a point
$I \in \mathcal{T}$, one defines an almost-complex structure
$\mathcal{I}$ on $T_I \mathcal{T}$ in the following way : the
connection $\nabla$ splits $T_I \mathcal{T}$ into the
space $T_I \mathcal{T}_{\pi(I)}$ tangent to the fibers and an
horizontal space 
$$Hor_I \mathcal{T} \simeq T_{\pi(I)} M  = H_{\pi(I)} \oplus
W_{\pi(I)}\, .$$
Changing the basis the local basis $(I_1, I_2 ,I_3)$, one can assume
that $I=I_1$, and the almost complex structure $\mathcal{I}$ is the
natural one on $T_{I_1}  \mathcal{T}_{\pi(I)} \simeq T S^2$,
whereas $\mathcal{I}|_H = I_1$ and $\mathcal{I} (R_2) = R_3 $. 

\begin{thm}[\cite{Biq00},\cite{Duc04}]
Let $(M^{4n+3}, H )$ be an integrable QC-manifold. The triple
$(\mathcal{T}, \eta , \mathcal{I} ) $ is a CR-integrable structure of
signature $(4n+2, 2 )$ and called the twistor space of $M$.
\end{thm}

One has the following result :

\begin{prop}[\cite{LeB82}, \cite{Biq00}]
Let $(M, H)$ be an analytic quaternionic contact distribution and let $(
\mathcal{T}, \eta, \mathcal{I})$ be its twistor space. There exists a
contact holomorphic manifold $ (\mathcal{N}^{2n+3}, \eta^c )$, a
family $(\mathcal{C}_m)_{m \in N^c } $ of dimension $2n+2$ of smooth rational
curves in $\mathcal{N}^c$ and a hypersurface $M^c \subset N^c$ such that :
\begin{itemize}
\item The distribution $\ker \eta^c $ is transverse  to the curves
$\mathcal{C}_m$ as soon as $m \in N^c - M^c$.
\item The normal bundle of the curves $\mathcal{C}_m$ is $\mathcal{O}
( 1) \otimes \mathbb{C}^{2n+2}$.
\item There exists a compatible real structure $\sigma $ on $\mathcal{N}$ such
that $M$ is the real slice of $M^c$ and $\mathcal{T}$ is a real
hypersurface in $\mathcal{N}$ with the induced CR-structure.
\item There exists an holomorphic metric $g$ on $N^c - M^c$ with holonomy $Sp_{n+1} (
\mathbb{C})Sp(1)$ and conformal infinity $H$. 
\end{itemize}
\end{prop}

\subsection{Embedding a QC structure in a quaternionic manifold}
Let $(M^{4n+3}, H )$ be an integrable QC-manifold. We use the
notations of section \ref{section-twist}. If $x\in N^c$, then $N_x$
stands for the normal
fiber bundle of $\mathcal{C}_x$. We put
$$ \begin{array}{rcl}
E_x & = & H^0 ( \mathcal{C}_x, L^{- \frac 12} \otimes N_x ) \, , \\
A_x & = & H^0 ( \mathcal{C}_x, L^{\frac 12} ) 
\end{array} 
$$
so that $T_x N^c = H^0 ( \mathcal{C}_x, N_x ) = E_x \otimes
A_x$. We get a $GL(2n+2, \mathbb{C} ) GL (2, \mathbb{C})
$-structure on $N^c$ and an almost-quaternionic structure on $N$ such that on
$N^c- M^c$, the $Sp_{n+1} ( \mathbb{C}) Sp(1)$-structure is a
reduction of this $GL(2n+2, \mathbb{C} )GL( 2
,\mathbb{C})$-structure. In particular, one sees that the almost
quaternionic structure of $N- M$ admits a smooth prolongation to
$M$. Because the almost-quaternionic on $N-M$ admits a
quaternionic-K\"ahler metric, its torsion vanishes on $N-M$, and so on
$N$ by continuity. 

If $x\in M^c$, the tangent space of the curve $\mathcal{C}_x$ lies in
the kernel of $\eta^c$. It follows that the hyperplane $H_x^c = H^0 (
\mathcal{C}_x, \ker \eta^c \,  / T \mathcal{C}_x ) \subset T_x M^c $ is
well defined. In fact, if $x \in M$, then $H_x$ is the real part of
$H^c_x$. In the decomposition $T_x N^c = E_x \otimes A_x$, we see that $H_x^c$
is the kernel of the 1-form $\eta = \eta^c \otimes 1 \in ( E \otimes A )^*
\otimes L $. We deduce that :

\begin{lm}
The hyperplane $H^c_x$ is invariant under the action of the subgroup
$GL(A)$ of $GL( T_x N^c)$.
\end{lm}

To  summarize, we have proved the following result:

\begin{thm}
Let $(M^{4n+3}, H) $ be an integrable QC-manifold. There exists a
quaternionic manifold $(N^{4n+4}, \mathcal{Q})$ such that :
\item[(i)] $M$ is a QC-hypersurface of $N$ and $H$ is the
$\mathcal{Q}$-invariant subspace of $TN$.
\item[(ii)] $M$ separates $N$ into two quaternionic manifolds $N_+$
and $N_-$ such that $N_+$ has a definite-positive
quaternionic-K\"ahler metric with conformal infinity $H$ and $N_-$ has
a quaternionic-K\"ahler metric with signature $(4, 4n)$ and conformal
infinity $H$. 
\end{thm}

We give now an explicit torsion-free connection on $N^c$ that preserves the
quaternionic structure. Let $\nabla$ be the Levi-Civita connection of
$g$. This is a meromorphic connection on $N^c$. Let $\rho$ be an holomorphic function defined on a neighbourhood of a
point $p \in N^c$ and vanishing up to order one on $M^c$. One knows that
we can write
$$g = \frac{1}{\rho^2}   ( (d \rho) ^2 + \eta_1^2 + \eta_2 ^2 +
\eta_3 ^2  ) + \frac{1}{ \rho} g_H  +\cdots  \, . $$

\begin{prop}
The connection $\nabla + a^{ \frac{d\rho}{\rho} } $ defined on a
neighbourhood of $p \in M^c$ is holomorphic and its restriction to $N$
gives a torsion-free connection preserving the quaternionic structure
of $N$.
\end{prop}

\begin{proof}
 We write
$$g = \frac{1}{\rho^2} ( ( d \rho)  ^2 + \eta_1^2 + \eta_2^2 +\eta_3^2
) + \frac{1}{\rho} g_{-1} + g_0 + \cdots $$
where $g_{-1} |H = g_H$ and $g_i$ is a covariant 2-tensor which does
not depend of $\rho$. Because $g_{-1}|_H$ is non degenerate, one can
define the holomorphic orthogonal $W$ of $H$ with respect to $g_{-1}$. The dots will indicate terms of order strictly
inferior in $\rho$ when $\rho$ goes to zero.

Let $(X_i)_{i \geq 4}  $ be an orthonormal basis of $H$ respectively
to $g_H$, $(X_i)_{1\leq i \leq 3} $ be an orthonormal basis of $W \cap
\ker d \rho$ for
$ \sum_i \eta_i^2$ and $X_0$ such that $d \rho ( X_0 ) = 1$,  $\eta_i
( \partial_{\rho} ) =0  $ and $g_{-1} ( X_0 , H) = 0$. The $Sp_{n+1} ( \mathbb{C} )Sp_1 ( \mathbb{C} ) $-structure is sent holomorphically to the
$Sp_{n+1} ( \mathbb{C} ) Sp_1 ( \mathbb{C} ) $-structure of the complexification of
$g_{\mathbb{H}}$, up to first order in $\rho$. Moreover, writing the
equalities of the form
$$g( [ X_i, X_j ] ,  X_k ) = - \frac{1}{\rho^2}\sum_p d \eta_p ( X_i,
X_j) X_k   + \cdots $$
if $ i, j \geq 4$ and $1 \leq k \leq 3$, one sees that it sends
the 1-jet of the $Sp_{n+1} ( \mathbb{C} )Sp_1 ( \mathbb{C} ) $ structure holomorphically to the
1-jet of the $Sp_{n+1} ( \mathbb{C} )$ $ Sp_1 ( \mathbb{C} ) $-structure of the complexification of
$g_{\mathbb{H}}$, up to first order in $\rho$.  Therefore $\nabla$ has
a pole of order 1 along $M^c$ and $\nabla + a^{ \frac{d
\rho}{\rho}}$ admits an holomorphic continuation to $N^c$.
\end{proof}

\subsection{Examples} In this section, we describe an illustration of
this theorem with the help of a family quaternionic-K\"ahler metrics obtained
by Bogdan Alexandrov in \cite{Bog01}. The construction begins with the data
of a hyper-K\"ahler manifold $ (M', I_1', I_2', I_3' , g) $ with
K\"ahler forms $w_i'$ satisfying the hypothesis that there exists 1-forms
$\alpha_i$ such that $d \alpha_i = w_i '$. Let $( \rho, x_1 ,x_2,
x_3)$ be the coordinates of $\mathbb{H}$, $N = \mathbb{H} \times M'$
and $M$ be the hypersurface defined by $\rho =0 $ and let $\pi : N
\rightarrow M'$ be the projection. One defines an hypercomplex structure $(I_1, I_2,
I_3)$ by the formula
\begin{equation}
I_i d\rho  = - d x_i + \pi^* \alpha_i \mbox{ and } I_i \pi^* \beta =
\pi^* I_i' \beta \mbox{ if } \beta \in T^* M' \, .
\end{equation}

Then, $M$ is a quaternionic contact hypersurface in $N$ with contact
distribution $H = \cap \ker I_i   d \rho $ and it is the
conformal infinity or the quaternionic-K\"ahler metric
$$ g = - \frac{1}{ \rho ^2} ( d\rho ^2 + \sum_i  ( I_i d\rho )^2 ) +
\frac 1 \rho  \pi^* g' \, .$$
Moreover, the compatible metric on $H$ is $g' \circ \pi_* | _H $.

\section{Second fundamental forms of a QC-hypersurface}

One may ask the following question: Does every embedding of a quaternionic-contact structure arises
in the way of Theorem 1.1 ? The aim of this section is to give a first
step toward an answer to this question. In particular, I will define
second fundamental forms for QC-hypersurfaces and prove that it
vanishes for the embedding given by Theorem 1.1.

\bigskip

Let $(N, \mathcal{Q}, \nabla )$ be a quaternionic manifold, $f$ a smooth
function on $N$ and assume that $M$ is a QC-contact
hypersurface with defining function $f$, QC-distribution $H$ and compatible definite positive
metric $g = 2 \nabla df |_H$. Let $(I_1, I_2, I_3)$ be a local choice of
quaternionic structure, $\xi$ be the normal vector field of $f$ along
$M$, $\eta_i = - df \circ I_i $ and $R_i = I_i \xi + r_i$ be the
Reeb vector field of $\eta_i$. From the proof of theorem 3.2, we
get that $r_i \in H$ is determined by $\nabla df ( I_i \xi + 2 r_i , X
) = 0$ for $X \in H$. The subvector bundle spanned by $I_1 \xi$, $I_2
\xi$ and $I_3 \xi$ is called $W^{\xi}$, and $p_H$ is the projection
onto $H$ with kernel $W^{\xi} \oplus \mathbb{R} \xi$.

\subsection{Covariant derivatives in the direction of $H$}

One has a natural connection on $TM$ defined by $\nabla^H _X Y =
p_H \nabla_X Y$ if $X, Y \in TM $. We compare now this connection with
the Biquard connection $\nabla^g$ on $H$.

The kernel of $p_H$ is
$\mathcal{Q}$-invariant, thus $\nabla^H$ preserves the $GL(n,
\mathbb{H} ) Sp(1)$-structure of $H$. Because the torsion of $\nabla$
is zero, for $X, Y \in H$, the torsion $T^H$ of $\nabla^H$ satisfies 
$$T^H (X, Y ) =  - \sum_i d \eta_i ( X, Y) I_i \xi \, .$$ 
Let $\mathfrak{a}$ the skew-symmetrisation and $X,Y$ and $Z$ be in $H$, so that if $ w_i = d
\eta_i |_H$, then 
$$ 0 = \mathfrak{a} ( \nabla^H w_i ) ( X, Y, Z) $$ 
$$ + d \eta_i( T^H ( X, Y)
, Z) + d \eta_i  (T^H ( Z, X) , Y ) + d \eta_i ( T^H ( Y, Z ) , X  )
\, , $$
and therefore if we restrict now ourselves on $H$, 
$$ \mathfrak{a} ( \nabla^H w_i ) = - \sum_j  i_{I_j \xi }
d  \eta_i |_H \wedge w _j \, . $$
Let $\Omega = \sum_i w_i^2$ be the fundamental form of the quaternionic
structure of $H$. We have obtained
\begin{equation} \label{eq antisymetrie}
\mathfrak{a} ( \nabla^H \Omega )  \in [ \sigma^2] \otimes [\sigma^2] \otimes [ \lambda^1
  \sigma^1 ] \cap \Lambda^5  = [ \lambda^1 \sigma^3]
  \oplus  [ \lambda^1 \sigma^1 ] \, . 
\end{equation}
The skew-symmetrisation $\mathfrak{a} : \Lambda^1 \otimes \Lambda^4
  \rightarrow \Lambda^5$ is injective if $n \geq 3$ and its kernel is
  $[\lambda^3_1 \sigma ^3 ]$ when $n=2$ so that $\nabla \Omega \in
  [\lambda^1 \sigma^3 ] \oplus [\lambda^1 \sigma^1 ] $ if $n \geq 3$
  and $\nabla \Omega \in [\lambda ^1 \sigma ^3 ] \oplus [ \lambda^1
  \sigma^1] \oplus [\lambda^3_1 \sigma^3]$ if $n=2$. On the other hand, $ \nabla^H$ is quaternionic, so that 
\begin{equation}
\nabla^H \Omega \in \Lambda^1 \otimes \left (  ( \mathfrak{sp}(n)
\oplus \mathfrak{sp}(1) )
^{\perp} \cap ( \mathfrak{gl}(n, \mathbb{H} ) \oplus \mathfrak{sp}(1)
)  \right )  
\,  ,
\end{equation}
i.e.
\begin{equation} \label{derivation}
\nabla^H \Omega \in  2 [ \lambda^1 \sigma^1 ] \oplus [ \lambda^3_1 \sigma^1] \oplus [ \lambda^3_0
\sigma^1 ] \, .
\end{equation}
Comparing (\ref{eq antisymetrie}) and (\ref{derivation}), we see that
$\nabla^H \Omega$ must live in a factor isomorphic to
$[\lambda^1 \sigma^1]$. A change of quaternionic connection $\nabla
\rightarrow \nabla + a^{\theta}$
changes $\nabla^H \Omega$ by a factor $\theta|_H \in [\lambda^1
\sigma^1]$. Hence, one can choose $\theta$ in such a way that $ \nabla^H
\Omega $ vanishes.

\begin{prop} \label{complementary}
 Assumes that $n \geq 2$. The adapted complementary vector bundle
 $W^g$ is equal to $W^{\xi}$. Moreover, one can choose $\nabla$ in such a way that $\nabla^H_X Y = \nabla^g_X Y$ for $X, Y \in H$.
\end{prop} 

\begin{cor}
Assumes that $n \geq 2$ and that $\nabla$ is chosen as in
 \ref{complementary}. The partial covariant derivatives $ \nabla^H$
 and $ \nabla^g :
 \Gamma ( \mathcal{Q} ) \rightarrow \Gamma ( H^* \otimes \mathcal{Q}
 )$ are equal. 
\end{cor}

\subsection{Covariante derivative in the direction of $W^g$}
In this subsection, we assume that $n \geq 2$ and that $\nabla$ is
chosen as in proposition \ref{complementary}. Let $\varepsilon^{ijk}$
be the signature of the permutation $(123) \rightarrow (ijk)$, and let
$(e_k)$ be an orthonormal basis of $H$, $(e_k^*)$ be its dual basis and $w_i  = \frac 12 \sum_k e_i^* \wedge  I_i e_k^* \in \Lambda^2 H^*$ be the
restriction $d \eta_i |_H$. We put $w_i^* = \frac 12 \sum_k e_i \wedge
I_i e_k$.

\begin{lm} \label{derivation dans la direction de W}
For all $i , j \in \{ 1, 2, 3  \}$, one has the following formula giving the action of $\nabla^g$ on
$\mathcal{Q}$ :
\begin{equation} 
\nabla_{I_i \xi  }^g I_j = \nabla_{I_i \xi }^H I_j + \sum_{p, k}
\varepsilon^{jpk} df ( R^{\nabla}_{w_k^* } I_i \xi ) ) I_p  - \sum_{p, k=1}^3
\varepsilon^{ jpk} \nabla df ( I_k \xi , I_i \xi ) I_p \, .   
\end{equation}
\end{lm}

\begin{proof}
Let $T$ be the torsion of the Biquard connection $\nabla^g$. If $X \in
H$, the
$H$-part of the torsion $T^H ( I_i \xi , X ) $ of $\nabla^H$ is equal
to $ \nabla ^H_{X} I_i \xi$. Therefore, we get :
$$( \nabla^g_{I_i \xi } I_j  - \nabla_{I_i \xi }^H  I_j  ) X = T ( I_i
\xi , I_j X ) - I_j T( I_i \xi , X ) - \nabla^H _{I_j X } I_i \xi +
I_j \nabla^H_{X} I_i \xi \, .    $$

Because both $\nabla^g$ and $\nabla^H$ preserve $\mathcal{Q}$, we can
now take the $\mathfrak{sp}(1)$-part of this expression. The
projection onto $\mathfrak{sp}(1)$ of $
 (T ( I_i
\xi , I_j \cdot ) - I_j T( I_i \xi , \cdot ))|_H$
 vanishes (see section 3.3), hence
$$  \nabla^g_{I_i \xi } I_j  - \nabla_{I_i \xi }^H  I_j  =  (-\nabla^H _{I_j \cdot} I_i \xi +
I_j \nabla^H_{\cdot } I_i \xi )_{\mathfrak{sp}(1)} \, .  $$ 

 Let $p \neq j$ be in $\{1 ,2 ,3 \}$. As
$\nabla df ( W^{\xi } \oplus \mathbb{R} \xi , H ) = 0 $, one
gets for $X \in H$,
$$\begin{array}{rcl}
 g (I_j \nabla^H_{X}I_i \xi - \nabla^H _{ I_j
X} I_i \xi  , I_p X ) &   = &   \nabla df  (I_j \nabla_{X}I_i \xi - \nabla _{ I_j
X} I_i \xi  , I_p X )  \\
& =  & - \nabla df (  I_j I_p X , \nabla_X  I_i \xi )\\
& &  - \nabla df (
 I_p X
 , \nabla_{I_j X } I_i \xi ) \, . \\

\end{array}    $$  

One has
$$\nabla df  (I_j  I_p X , \nabla_X  I_i \xi  ) = ( I_j I_p X ) . df
 ( \nabla_X  I_i \xi  ) - df ( \nabla_{I_j I_p X} \nabla_X  I_i \xi
 ) = $$ $$
 (I_j I_p  X) . ( - \nabla df ( X,  I_i \xi ) + X . df ( I_i \xi
 ) )  - df ( \nabla_{I_j I_p X} \nabla_X  I_i \xi
 )  $$
$$ = - df ( \nabla_{I_j I_p  X } \nabla_X I_i \xi ) \, , $$
so that
$$ \begin{array}{rcl}
 g (I_j \nabla^H_{X}I_i \xi - \nabla^H _{ I_j
X} I_i \xi  , I_p X ) &   = &  df ( \nabla_{I_j I_p X} \nabla_X I_i \xi +
\nabla_{I_p X} \nabla_{I_j X } I_i \xi ) \, . \\

\end{array}$$
We replace now $X$ by $e_i$ and sum over the basis $(e_1,\dots, e_{4n} )$ of $H$ to obtain the lemma.
\end{proof}

\subsection{Second fundamental forms} By a conformal change $ \eta_i
'=  g^2 \eta_i$, the tensor $w_i$ becomes $( w_i^* )' = g^{-2} w_i^*$.
\begin{defi}
The conformal second fondamental form of the embedding $M \hookrightarrow N$ is the
trace-free part of the tensor
$$\sum_{i ,j} ( \nabla df ( I_i \xi , I_j \xi  )- 4 \nabla df (r_i ,
r_j ) ) \otimes w_i \otimes I_j \, . $$
It does not depend of $f$ nor of the torsion-free connection $\nabla$.
\end{defi}

\begin{rem}
When $n=0$, we obtain the trace-free part of the second fundamental
form of the embedding of a manifold $M$ into a conformal manifold $(N
, [g] )$. 
\end{rem}

\begin{proof}
Assume that $f' = g^2 f$ is another defining function of $M$ where $g$
does not vanish, and that $\xi '$ is the normal vector field of $f'$
along $M$. Then, one has $df'
= g^2 df $ and $\nabla df'  = 2 g dg \otimes df + 2g df \otimes dg +
g^2 \nabla df $ along $M$. Hence, if $\xi ' = g^{-2} \xi +r  $ with $d
f'( r) = 0$, then the formula $ d f' ( I_i \xi' ) =  df ( I_i \xi ) +
g^2 df ( I_i r) = 0 $
implies that $r \in H$. The Reeb vector fields $R_i'$ of the forms
$\eta_i' = g^2 \eta_i$ can be written $R_i ' = I_i \xi ' + r_i '  $
and satisfy  $ \nabla d f' (I_i \xi ' + 2 r_i ', X   )= 0    $ for $X \in
H$, whence
$$r_i  ' = g^{-2} r_ i - \frac 12 I_i r \, . $$ 
We obtain
$$\begin{array}{rcl}
 \nabla df' ( I_i \xi'  , I_i \xi'   )+ 4 \nabla df' (r_i ,
r_j )  & =  & g^{-2} ( \nabla df ( I_i \xi , I_j \xi ) -4 \nabla df ( r_i , r_j
 )   )  \\
& & +  \nabla df ( I_i \xi
, I_j r )  + \nabla df ( I_i r , I_j \xi  ) \\
& &  + 2 \nabla df ( r_i , I_j r
 ) + 2 \nabla df ( r_j , I_i r ) \, .
\end{array}
 $$
Therefore, the conformal curvature is independant of $f$. A quick
computation and the formula of lemma 2.1 give that it is independant of the connection $\nabla$.

\end{proof}

\begin{defi} Assume $n=1$. The projection $Q$ onto $[\lambda^1
    \sigma^3]$ of the tensor 
$$\sum_i \nabla df ( I_i \xi, \cdot) |_H \otimes I_i
\in [\lambda^1_0 \sigma^1 ] \otimes [ \sigma ^2 ] \simeq [ \lambda^1_0
\sigma ^1 ] \oplus [ \lambda^1_0 \sigma^3  ]$$
 is called the horizontal second fundamental form of the
QC-hypersurface $M$ in $N$. It does not depend of the choice of
torsion free connection $\nabla$ preserving $\mathcal{Q}$ nor of the
choice of defining function $f$ for $M$.
\end{defi}  

\begin{proof}
Using the proof of the definition 3.1, we obtain that for $X \in H$,
 $$\nabla df' ( I_i \xi ', X) = \nabla df ( I_i \xi,
X)   +g^2 \nabla df (  I_i r  ,X ) $$
and finally that only the factor in $[ \lambda^1_0 \sigma^1 ] $ changes.
\end{proof}

\subsection{Second fundamental forms and AHQK-metrics}
\begin{thm}
Assume that $(M, H)$ is a QC-hypersurface in a quaternionic manifold
$(N, \mathcal{Q})$ such that there exists a AHQK-metric $g$ defined on
an open subset of $N$, compatible with $\mathcal{Q}$ and with conformal infinity $(M,
H)$. Then the second fundamental forms of $M$ vanish.
\end{thm}

\begin{proof}
 Let
  $\rho$ be a defining function of $M$, let $\nabla$ be a
  quaternionic connection such that $\tilde{\nabla}= \nabla + a^{ - \frac{d \rho
      }{\rho } } $ is the Levi-Civita connection of $g$. Let $(I_1,
  I_2, I_3)$ be a local trivialization of $\mathcal{Q}$ around a point
  $p \in M$ such that the almost complex structures $I_i$ are parallel
  at $p$ respectively to $\nabla$ and put $\nabla I_i = \sum_{j,k}
  \varepsilon^{ijk} \beta_{k} \otimes  I_j$.

We compare the curvature $R^{\tilde{\nabla}}$ and $R^{\nabla}$ acting
  on $\mathcal{Q}$. We have the well known formula
$$R^{\tilde{\nabla}} = R^{\nabla} + d^{\nabla} a^{- \frac{d
      \rho}{\rho} } + [ a^{- \frac{ d \rho } {\rho} }  , a^{- \frac{d
      \rho}{\rho }} ] \, , $$
which gives
$$R^{\tilde{\nabla}}  - R^ {\nabla}  =   \sum_{i=1}^3  \frac{1}{\rho^2} \left ( d (\rho \circ I_i ) \wedge
d \rho  + \frac 12 \sum_{j, k} \varepsilon^{ijk} d \rho \circ I_j
    \wedge d \rho \circ I_k  \right ) \otimes I_i $$
$$ + \sum_{i=1}^3   \frac{1}{\rho} d ( d \rho
    \circ  I_i )  \otimes I_i  \, . $$
Let $w_i$ be the K\"ahler form $w_i (\cdot ,\cdot ) = g ( I_i \cdot ,
\cdot ) $. The metric $g$ is quaternionic K\"ahler, hence
$R^{\tilde{\nabla}} = \sum_i w_i \otimes I_i$ on $\mathcal{Q}$. The
term $R^{\nabla}$ in the previous formula extends smoothly on $M$, so
that we obtain
$$w_i =  \frac{1}{\rho^2} \left ( d (\rho \circ I_i ) \wedge
d \rho ) + \frac 12 \sum_{j, k} \varepsilon^{ijk} d \rho \circ I_j
    \wedge d \rho \circ I_k  \right ) + \frac{1}{\rho} d ( d \rho
    \circ  I_i ) + \gamma_i $$
where $\gamma_i$ extends smoothly on $M$. In
particular, we get for all $i$,
$$g   = \frac{1}{\rho^2} ( d \rho + \sum_{j=1}^3 ( I_j d \rho)^2
)   $$ $$
 -
\frac{1}{\rho } ( \nabla d \rho  + \nabla d \rho ( I_i  \cdot , I_i
\cdot 
)
  + \sum_{j,k} \varepsilon^{ijk} ( \beta_k \otimes d \rho  +  I_i d
\rho  \otimes I_i \beta_k )  ) + \cdots   $$
where the dots indicates terms that extend smoothly on $M$. We deduce
that at $p$, one has $\nabla d \rho ( X, Y ) = \nabla d \rho ( I_i X, I_i Y ) $
for all $i$ and $X$, $Y$. Let $\xi$ be the normal vector field of
$\rho$ along $M$. It turns out that it is equivalent to
$W^{\xi} = W^g$ and to the vanishing of the second fundamental forms.
\end{proof}

\section{Final remarks}

This paper is a first step in the direction of understanding the
structures arising on hypersurfaces in quaternionic manifolds. In
the case of quaternionic contact structures, it is interesting
to know if one can construct an embedding of a QC-distribution into a
quaternionic manifold with a given second fundamental form in the same
way that in the
$3$-dimensional case (\cite{LeB85}).  We suspect that there are obstructions involving higher derivatives of the second fundamental forms,
obstructions that one should be able to recognize in the construction of an
adapted twistor space.

In the general case of weakly QC-distribution, a CR-twistor space
should exist, but the construction of an adapted Biquard connection
remains to be done. We will adress these problems in a future work.

\end{document}